\documentclass[11pt]{amsart}
\usepackage{amscd,amsmath,amssymb,amsthm}
\usepackage{color}
\usepackage{textcomp}
\usepackage{graphicx}
\usepackage{tikz}

\newtheorem{theorem}{{\bfseries Theorem}}[section]

\newtheorem{proposition}[theorem]{{ Proposition}}
\newtheorem{lemma}[theorem]{{ Lemma}}
\newtheorem{corollary}[theorem]{{ Corollary}}
\newtheorem{definition}[theorem]{{ Definition}}
\newtheorem{example}[theorem]{{ Example}}
\newtheorem{remark}[theorem]{{ Remark}}

\newcommand{\R}{\mathbb{R}}

\newcommand{\N}{\mathbb{N}}
\newcommand{\Z}{\mathbb{Z}}

\newcommand{\Q}{\mathbb{Q}}

\begin{document}

\title{On Almost periodicity and minimality for semiflows}

\author{Joseph Auslander}
\address{Department of Mathematics, University of Maryland,  College Park, MD 20742, USA.}

\email{ jna@umd.edu}

\author{Anima Nagar}
\address{Department of Mathematics, Indian Institute of Technology Delhi,
	Hauz Khas, New Delhi 110016, INDIA.}

\email{ anima@maths.iitd.ac.in}

\thanks{We thank a very generous, anonymous reviewer for all the comments and suggestions that have enhanced the write-up of this article.}

%\date{20 December, 2020}

\vspace{0.5cm}


\maketitle

\begin{abstract}
In topological dynamics, the dynamical behavior sometimes has a sharp contrast when the action is by semigroups or monoids to when the action is by groups. In this article we bring out this contrast while discussing  the equivalence of  almost periodicity and minimality, and some implications when every point is an almost periodic point.
	
\end{abstract}

{\textsc{keywords}:}  flows, semiflows, minimality, almost periodic points.

\vskip 1cm

Let $X$ be a compact, Hausdorff topological space and $T$ be a topological group or  monoid (semigroup with identity).  We say that $T$ acts on $X$ when there is a continuous map $\pi: T \times X \to X$ such that  $\pi (e, x) = x$ for all $x\in X$, where $e \in T$ is the identity in the monoid or group $T$;
satisfying the Kolmorogov identity: $\pi (ts, x) = \pi (t, \pi (s,x))$ for all $t, s \in T$ and all $x\in X$.

\bigskip

$(X,T,\pi)$ or $(X,T)$ is called a \emph{topological dynamical system ( TDS )}.

\medskip

Here $X$ is called the \emph{phase space}, $T$ the \emph{acting group or monoid} and the \emph{action} $\pi$ gives the homeomorphism or continuous map $\pi^t: X \to X$ and we write    $ tx$ for $\pi^t(x)$.

When $T$ is a topological group - $(X,T)$ is called a \emph{flow}, and for a topological \emph{monoid or semigroup} $T$ - $(X,T)$ is called a semiflow. 

When the acting group $T = \Z$ or the monoid $T = \Z_+$ or the semigroup $T = \N$, $\pi^1 = f$ gives a generating homeomorphism or continuous map on $X$, i.e. $f(x) = \pi(1,x)$, which gives iterations $f^n(x) = \pi(n,x) = nx$. In this case we call the system $(X,f)$ a \emph{cascade} or a \emph{semicascade} respectively.

 \bigskip
 
 Let   $x \in X$, then $Tx = \{tx: t \in T\}$ is called the \emph{orbit}, and $\overline{Tx}$ denotes the \emph{orbit closure} of $x$. For a cascade or semicascade $(X,f)$, the \emph{(forward)orbit} of the point $x$  is $O(x) = \{f^n(x) : n \in \N\}$.
 
 \bigskip
 
 For a $ TDS $ $(X,T) $, a  subset $A  \subseteq X$ is called \emph{invariant} if $TA \subseteq A$. Note  that for any  $A  \subseteq X$, 
  the subset $ TA $ is invariant and that the closure of an invariant set
  is invariant. And  a subset $M  \subseteq X$ is called \emph{minimal} if it is non-empty, closed, invariant and minimal with respect to these properties. Thus $M$ is minimal if and only if $\overline{Tx} = M$ for every $x \in M$. And Zorn's lemma guarantees that every $ TDS $ admits a minimal subset.

\bigskip

Fundamental to the study of the topological dynamics of group actions is the equivalence of almost periodicity and minimal sets. 

\bigskip

For any group or monoid or semigroup $T$, a subset $A \subset T$ is called \emph{syndetic}  if there exists a compact $K \subset T$ such that $T = K \cdot A$ (see \cite{gh}). Also,  a semiflow $(X,T)$  is called \emph{surjective} if each $t \in T$ is surjective.

\bigskip

Let $(X,T)$ be a flow. For $ x \in X $ and $ U $ a neighbourhood of $ x $, define $ A_U(x) :=\{t \in T: tx \in U\} $. We just mention $A_U$ whenever $U \ni x$ is understood.

A point $ x \in X$ is said to be an \emph{almost periodic point} for the flow $(X,T)$ if for every neighbourhood $ U \ni x $, there is a compact  $ K \subset T $ with $ T=K \cdot A_U $, i.e.  for every open $U \ni x$ the set $A_U$ is syndetic.

\bigskip

We review the proof of the equivalence of almost periodicity  and minimality of the orbit closure for the flow $(X,T)$.

\begin{theorem} For a TDS $(X,T)$ with $T$ a group and $x \in X$, the following are equivalent.
	\begin{itemize}
		\item[(i)]   The orbit closure $\overline{Tx}$ is a minimal subset.
		
		\item[(ii)] For every open $U$ with $x \in U$, there exists a compact subset $K$ of $T$ such that $T = K \cdot A_U$, i.e. each $A_U$
		is syndetic.
		
		\item[(iii)] For every open $U$ with $x \in U$, there exists a finite subset $F$ of $T$ such that $T = F \cdot A_U$,  i.e. each $A_U$
		is syndetic when $T$ is given the discrete topology.
\end{itemize} \end{theorem}

\begin{proof}
	
We note that $(iii) \implies (ii)$ is a fortiori. 

We now prove that $(ii) \implies (i)$. Let $ x $ be an almost periodic point and let $ M $ be a minimal set contained in $ \overline{Tx} $. Let $ U $ be a neighbourhood of $ x $, so $ T=K \cdot A_U $ with $ K $ compact. Let $ \{t_i\} $ in $ T $ with $ t_ix \to m \in M $. Then $ t_i=k_ia_i $ with $ k_i \in K $ and $ a_i \in A_U $. We may suppose $ k_i \to k \in K $ and $ a_ix \to x' \in \overline{ U}  $  so that $ m=kx' $ and $ x'=k^{-1}m \in M $. Since U is an arbitrary neighborhood of $ x $ and $ M $ is closed it follows that $ x \in M $ and $ \overline{Tx}=M $.

 As for $(i) \implies (iii)$, we note that this is a classical result and    it will be shown to follow from a semiflow theorem below (Theorem \ref{msf}).

\end{proof} 

\bigskip

 In this note we show that if the orbit closure of $ x $ is minimal for a semiflow, it need not be the case that $ x $ is almost periodic. However with a modified definition of almost periodicity (which coincides with the above definition in the case of an acting group) we can recover the equivalence of  almost periodicity and minimality. We further discuss more implications when every point is an almost periodic point.

\section{Some deviations}

 When $T$ is a group, a point $x \in X$ is an almost periodic point if and only if $\overline{Tx}$ is a minimal subset of $X$.  This further guarantees that the flow $(X,T)$ admits an almost periodic point.

 \bigskip
 
 We present  examples which show that for semiflows $(X,T)$ neither of the properties of almost periodicity  and minimality of
 the orbit closure implies the other - it may not be true that all points that are almost periodic points (satisfying $T = K \cdot A_U$) give minimal orbit closures and  minimal orbit closures need not imply that its points are almost periodic.

 \begin{example} \label{napp}
 	Let $f_\alpha: [0,1] \to [0,1]$, for $ \alpha \in [0,1]$ be defined as $f_\alpha(x)= \alpha x$.
 	
 	Consider $T = \{f_\alpha: \alpha \in [0, \frac{1}{2}] \cup \{1\}\}$. Then $T$ is a monoid with the uniform topology, under the operation of composition, acting on $[0,1]$ giving a semiflow $([0,1], T)$.
 	
 	We note that $T$ is a compact subset of the set of all continuous mappings on $[0,1]$ with the uniform topology.
 	
 	\bigskip
 	
 	Let $x = \frac{1}{2}$ and consider the open interval $U = (\frac{1}{3}, \frac{2}{3})$. Recall that $A_U = \{f_t \in T: f_t(x) \in U\}$ and note that here $A_U = \{f_1\}$, the singleton identity. Also $T = T \cdot A_U$, and so $A_U$ is syndetic and $x = \frac{1}{2}$ is almost periodic.
 	
 	However, $\overline{Tx} = \{\frac{s}{2}: s \in [0, \frac{1}{2}] \cup \{1\}\}$ is not minimal since $\{0\}$ is a minimal subset of $ \overline{Tx} $.
 \end{example}

In fact, we can say more on this:

\bigskip

For $t$ an element of a semigroup $T$ let $\rho_t$ and $\ell_t$ be the right and left translations given by
$\rho_t(s) = st,$ and $\ell_t(s) = ts$. 

For a semigroup $T$ the translation action of $T$ on itself is given by $\pi(t,s) = ts$.

\begin{theorem} \label{akin} For a monoid  $T$, the following conditions are equivalent:
	\begin{itemize}
		\item[(i)] For every $t \in T$, $\rho_t$ is surjective, i.e. $Tt = T$.
		
		\item[(ii)] For every $t \in T$, $\ell_t$ is surjective, i.e. $tT = T$.
		
		\item[(iii)] $T$ is a group.
	\end{itemize}
	Note that if $T$ is compact, (i) says that the translation action is minimal, while (ii) says that
	the translation action is surjective.
	
	If $T$ is a group, then both of the following hold.
	\begin{itemize}
		\item[(i)] For every $t \in T$, $\rho_t$ is injective.
		
		\item[(ii)] For every $t \in T$, $\ell_t$ is injective.
	\end{itemize}
	If $T$ is compact, then either of these conditions implies that $T$ is a group.
\end{theorem}

\begin{proof} For a group every $\rho_t$ and $\ell_t$ is bijective.
	
	If every $\rho_t$ is surjective, then  every  $t$ has a left inverse $\bar t$ such that $\bar t t = e$ and a $\tilde t$ such that
	$\tilde t \bar t = e$.  Hence, $\tilde t = \tilde t \bar t t = t$ and so $\bar t$ is the inverse of $t$.  For $\ell_t$ use the opposite
	monoid operation on the set $T$.
	
	If $u$ is an idempotent in $T$, then $u = u u = e u = u e$ and so if $\rho_u$ or $\ell_u$ is injective then $u = e$.
	
	If $T$ is compact, then every ideal $Tt$ contains an idempotent and so either (i) or (ii) implies that $e \in Tt$ and so
	every $t$ has a left inverse.  As above this implies that $T$ is a group.
	
\end{proof}

\begin{remark}Note that if  $T$ be a  monoid, then $\{ e \}$ is a syndetic subset if and only if $T$ is compact. 	
	
	%This obvious since $T = T \cdot \{ e \}$.
	
	Also when $T$ is a compact monoid and $(T,T)$  the TDS with the translation action,  then every point $t \in T$ is almost periodic.
	On the other hand, the action is minimal if and only if $T$ is a group.

\end{remark}

Thus, if $T$ is a compact monoid which is not a group then the translation action is non minimal but has every point almost periodic. As is illustrated in Example \ref{napp}, where one can use $T = [0,1]$ under multiplication.

\bigskip

\begin{example} \label{constant} Let $X$ be a compact metric space with no isolated points. 
	We consider $f_x (y) = x$ on $X$, which gives the constant map for every $x \in X$.  This makes $X$ a non-commutative semigroup with respect to the operation of composition of maps on $X$ giving the operation $yx = y$ on $X$. Let  $X_+ = X \cup \{ e \}$ where $e$ acts as an identity on this semigroup $X$.
	
	Let $D \subset X$ be a countable, dense subset. And consider the monoid $D_+ = D \cup \{e\}$.
	
	The TDS $(X,D_+)$ is minimal. However, for every $x \in X$, and open $U \ni x$, we have $A_U = (U \cap D) \cup \{e\}$. 
	
	Note that $D_+ = K \cdot A_U$ if and only if $K = D_+$. But $D_+$ is not compact (not closed). Thus, no $x \in X$ can be almost periodic.
	
\end{example}

\bigskip 

\begin{remark} Note in addition that although $(X,D_+)$ is minimal in the example above, we have $D_+ \cdot X = X $ and $  D \cdot X = D$ and so the action is far from surjective. 

In this regard we  note that, if $(X,T)$ is a minimal TDS, with $T$ a commutative monoid, then the action is surjective. This easily follows since for any $t \in T$, the subset $tX$ is closed invariant because $T$ is commutative and so $tX = X$.
	
\end{remark}

Recall that for a semicascade $(X,f)$, the \emph{omega limit set} of $x$ is defined as $\omega(x) = \{y \in X: $ there exists a seq $\{n_k\} \nearrow \infty$  such that $f^{n_k} x \to y\}$, and is the set of all limit points of  $\mathcal{O}(x)$. The point $x \in X$ is called \emph{recurrent} if $x \in \omega(x)$. 

\begin{lemma} \label{rec}
	For a semi-cascade $ (X,f) $, if every point in $X$ is recurrent then $ f $ is surjective.
\end{lemma}
\begin{proof}
	Let $ x \in X $ , and $\{ n_i\} \nearrow \infty $ with $ f^{n_i}(x) \to x $. A subsequence $ f^{n_i -1}(x) \to x' $. Then $ f(x')=x $.
\end{proof}

\bigskip
 
 We consider another surjective example with a non commutative monoid:
 
 \begin{example} \label{ex}
 	
	Consider  $X = \mathbf{T}$,	the unit circle bijective with $[0,1)$ and consider
 	$$T = \left\lbrace f_{n,q}(x) := nx+q \ (\mod 1): n \in \N \setminus \{1\}, \  q \in \mathbf{Q} \cap [0,1)\right\rbrace  \cup \{e\}$$
 	
 	with operation as composition of maps where $e$ is the identity map. 
 	
 	As can be seen $f_{m,r} \cdot f_{n,q} =  \begin{cases}
 	 f_{mn,mq+r}, \ mq+r < 1\\  	 
 	 f_{mn,mq-1+r}, \ mq+r \geq 1. \end{cases}$
  
   and so $T$ is a non commutative  monoid endowed with the uniform topology.
 	
 	And  the monoid $ T $ acts surjectively on $ X $ and so $ (X,T) $ is a surjective semiflow.
 	
 	\medskip
 	
% 	 We note that  for any $a \in [0,1) \cap \mathbf{Q}$, we have $a = 2a + (1-a) (\mod 1)$. Thus for the semiflow $(X,T)$, for each rational $y \in X$, we have $y  = f_{2,1-y}(y) = e(y)$. 

 For all $x \in X$, we see that
 	$ {Tx} \supset \bigcup \limits_{q \in \Q \cap[0,1)}\{f_{2,q}(x)\}$
 	 which is dense in $X$ and so $\overline{Tx} = X$. Thus $(X,T)$ is minimal.
 	
 	\medskip
 	
 	For every $n \in \N \setminus \{1\}$, let $T_n = \{f_{n,q}: q \in [0,1)\}$. Then
 	$ T = [\bigcup \limits_{n=2}^\infty T_n] \cup \{e\}$
 	is not equicontinuous.
 	
 	\bigskip
 	
 	Consider $0 \in X$ and let open $U \ni 0$ be an $\epsilon-$ neighbourhood for a small $\epsilon > 0$. Note that for all $n \in \N \setminus \{1\}$, $f_{n,0}(0) = 0$. Then $A_U =  \{f_{n,q} \in T: n \in \N \setminus \{1\}, q \in [0,\epsilon) \cup (1 - \epsilon,1) \} \cup \{e\}$. 
 	
 	\bigskip
 	
 	 	Suppose that  $A_U$ is syndetic. Then there  exists compact $K \subset T$ such that  $T = K \cdot A_U$.

 	\medskip
 	
 	Now for every $p \in [0,1)$, and $ f_{2,p} \in T$ there exists $ k_p \in K $ and $a_p \in A_U$ such that $ f_{2,p} =k_p \cdot a_p$. Then either $k_p=e$ or $a_p=e$  else $f_{2,p} = f_{mn,p} = f_{m,r} \cdot f_{n,s}$ with both $m, n \geq 2$, which is impossible. For $p \in (\epsilon, 1 - \epsilon)$ we note that $k_p \neq e$ else $t = a_p$ which contradicts that $f_{2,q} \in A_U$ only for $ q \in [0,\epsilon) \cup (1 - \epsilon,1)$. Hence $a_p = e$ for all  $p \in (\epsilon, 1 - \epsilon)$.
 	
 	Thus, for all $p \in (\epsilon, 1 - \epsilon)$ we have $f_{2,p} = k_p \in K$. Now let $r_n \in (\epsilon, 1- \epsilon) \cap \Q$, for $n \in \N$, be such that $r_n \to i$ with $i \in (\epsilon, 1-\epsilon) \cap (\R \setminus \Q)$. Define $g(x) = 2x +i \ (\mod 1)$ on $X$, then $f_{2,r_n} \to g$ pointwise. Now since the set $\{f_{2,r_n}: n \in \N\} \cup \{g\}$ is equicontinuous, we have $f_{2,r_n} \to g$ uniformly implying that $g \in K$ since $K$ is compact. But $g \notin T \supset K$ giving a contradiction. 
 	
  	Hence  $ A_U$ cannot be syndetic.

 	 	 \medskip
 	 
 	  Thus $0 $ is not an almost periodic point with the usual syndetic definition of almost periodic points given for flows.

 \end{example}

\bigskip  

\section{Minimal Semiflows}

We now give a modified definition of almost periodicity which in fact is equivalent to minimality in the case of semiflows.

\begin{definition} \label{got}
	For a  semiflow $(X,T)$, $x \in X$ is called a \emph{semiflow almost periodic point} if for every open $U \ni x$ and the set $A_U := \{t\in T: tx \in U\}$, there exists a finite $K \subset T$ such that $Kt \cap A_U \neq \emptyset$ for all $t \in T$.
\end{definition}

\begin{remark}
	We note that with this definition, all points in the Examples \ref{constant} and \ref{ex} from the previous section are semiflow almost periodic points. 
	
	For Example \ref{constant} we just consider a compact $K \subset D \cap U$, then $Kz \cap A_U \neq \emptyset$ for every  $z \in D_+$ there.
	
	Whereas for Example \ref{ex} we take $N \in \N$ large enough so that $\frac{1}{N} < \epsilon$. Again noting that for any rational $y \in X$ we have $f_{2,1-y}(y) =y$, we take $K = \{f_{2,0}, f_{2,\frac{1}{2N}}, \ldots, f_{2,\frac{2N-1}{2N}}\}$ so that for any $x \in X$, $Kx$ is $\frac{1}{2N}-$dense in $X$. Then for every $t \in T$, there exists $k \in K$ such that $kt \in A_U$, i.e. $kt(0) \in U$.
\end{remark}

\begin{remark}
	
	Define $k^{-1}(A)$ to be the preimage $\ell_k^{-1}(A)$ so that $K^{-1}(A) = \{ t : Kt \cap A \not= \emptyset \}$. Notice that
	if $T$ is a group then $K^{-1}(A) = K^{-1} \cdot A$.
	
	\medskip
	
	We note that when $T$ is a group, i.e., for the flow $(X,T)$, the condition $Kt \cap A \neq \emptyset$ reduces to $t \in K^{-1}A $ for all $t \in T$, which implies $T=K^{-1} \cdot A$, making the set $A$ syndetic.
	
	Conversely, if $ T = K \cdot A $ then $K^{-1}t \cap A \neq \emptyset$ for all $t \in T$ and thus the two definitions of almost periodic points are equivalent for flows.
	
	\bigskip
	
	This can be summarized as:
	
	\begin{displaymath}\begin{split}
			t \in s A_U \quad \Leftrightarrow \quad \exists \quad a \in A_U, t = sa \quad \Rightarrow \quad tx = sax \quad \Rightarrow \quad tx \in s U. \\
			\text{ if } \ \ T \ \ \text{is a group, then} \quad  tx \in s U \quad \Rightarrow \quad s^{-1}t \in A_U \quad \Rightarrow \quad t \in s A_U.
	\end{split}\end{displaymath}\vspace{.25cm}
	
	So if $K \subset T$ is compact, and $V$ is open containing $\overline{U}$, then
	\begin{displaymath}\begin{split}
			T = K A_U  \quad \Rightarrow \quad Tx \subset K \cdot U  \quad \Rightarrow \quad \overline{Tx} \subset K \cdot \overline{U}.\\
			\text{ If } \ \ T \ \ \text{is a group, then} \ Tx \subset K \cdot U \quad \Rightarrow \quad T = K A_U. \\
			\text{ If } \ \ T \ \ \text{is a group, then} \ \overline{Tx} \subset K \cdot \overline{U} \quad \Rightarrow \quad \overline{Tx} \subset F \cdot V, \\
			\text{ for some finite } F \subset K.
	\end{split}\end{displaymath}\vspace{.25cm}
	
	\begin{displaymath}\begin{split}
			t \in s^{-1}(A_U) \quad \Leftrightarrow \quad st \in A_U \quad \Leftrightarrow \quad stx \in U \quad \Leftrightarrow \quad tx \in s^{-1}(U). \\
			T = K^{-1}(A_U) \quad \Leftrightarrow \quad Tx \subset K^{-1}(U). \hspace{3cm} \\
			\text{If} \ \ F \subset T \ \ \text{is finite, then } \ \ Tx \subset F^{-1}(U) \quad \Rightarrow \quad \overline{Tx} \subset F^{-1}(\overline{U}).\\
			\overline{Tx} \subset K^{-1}(V) \quad \Rightarrow \quad \overline{Tx} \subset F^{-1}(V) \ \  \text{ for some finite } F \subset K.
	\end{split}\end{displaymath}
	
\end{remark}

\begin{theorem} \label{msf}
	For a semiflow $(X,T)$ and $x_0 \in X$, $ \overline{Tx_0} $ is minimal if and only if   $x_0$  is semiflow almost periodic.
\end{theorem}

\begin{proof}
Let $(X,T)$ be minimal and $x_0 \in X$. For open $ U \ni x_0 $ consider the corresponding $ A_U $.

For every $x \in X$, there exists $\sigma_x \in T$ and open $V_{\sigma_x} \ni x$ such that $\sigma_x V_{\sigma_x} \subset U$. By compactness, $ X=V_{\sigma_1} \cup \ldots \cup V_{\sigma_n} $.
Let $ K=\{\sigma_1,...,\sigma_n\} $ and let $ t \in T $. Then $ tx_0  \in$ some $ V_{\sigma_j} $, and so $ \sigma_j tx_0 \in U $. That is $ \sigma_j t \in A_U $. Thus $Kt \ \cap A_U \neq \emptyset$, and so $x_0$ is semiflow almost periodic.

Conversely, let $x_0 \in X$ be semiflow almost periodic and an open neighborhood $U  \ni x_0$. Let $U  \supset \overline{V}  \supset V \ni x_0$, and  $  K $ be a finite subset of $ T $ such that  $Kt \cap A_V \neq \emptyset$ for every $ t \in T $. 

Let $ y \in  \overline{Tx_0}$ and let $ s_j \in T $ with $ s_jx_0 \to y $. So for every $ s_j $ there exists a $ k_j $ such that $  k_js_j \in A_V $ and we may assume $ k_j \to k $ in $ K $.  Then there is  a subnet of $ ks_j \in A_V $ and so $ ks_jx_0 \in V $ and $ ky \in  \overline{V} \subset U $. Thus $ x_0 \in \overline{Ty} $ which proves the minimality of $ \overline{Tx_0} $.

\end{proof}

This improves Theorem 1 in \cite{go} (with the improved definition). Though for all practical purposes one can work with a compact $K$ in Definition \ref{got} as in \cite{go}.

\begin{corollary} For a TDS $(X,T)$ with  $x \in X$, the following are equivalent.
	\begin{itemize}
		\item[(i)]   The orbit closure $\overline{Tx}$ is a minimal subset.
		
		\item[(ii)] For every open $U$ with $x \in U$, there exists a compact subset $K$ of $T$ such that $T = K^{-1}(A_U)$.
		
		\item[(iii)] For every open $U$ with $x \in U$, there exists a finite subset $F$ of $T$ such that $T = F^{-1}(A_U)$.
\end{itemize} \end{corollary}

\begin{corollary} \label{admit}
	Every semiflow $(X,T)$ admits  semiflow almost periodic points.
\end{corollary}
\begin{proof}
	This follows since every semiflow $(X,T)$ admits minimal sets and each point in such a set is semiflow almost periodic.
\end{proof}

\begin{remark}
	From these results, we can replace the original topology on $T$ by the discrete topology.
\end{remark}

We show that for semicascades, both the definitions
of almost periodicity coincide, so in fact an orbit closure is minimal if and only if
all of its points are almost periodic in either sense. 

\begin{lemma} \label{equal}
	Let $ A \subset \Z_+ $ with $ 0 \in A  $ and let $ R>0 $. Then the following are equivalent:
	
	(i)  $ \Z_+=\{0,1, \ldots,R\}+A $.
	
	(ii) If $ m \in \Z_+ $ then $ (m+\{0,1,\ldots,R\}) \cap A \neq \emptyset$.	
\end{lemma}

\begin{proof}
	Suppose $  \Z_+ = \{0,1,2,..,R\} + A $ then for any $m \in \Z_+$, $ m+R=k+a\ $ with $ k \in \{0,1,2,...,R\} $ and $ a \in A $. So $  m+R-k=a \in A  $, i.e. $ m +  \{0,1,....,R\}  \cap A \neq \emptyset$ for all $ m \in \Z_+ $.
	
	Conversely, let $ m +  \{0,1,....,R\}  \cap A \neq \emptyset$ for all $ m \in \Z_+ $.  If $ m \in \{0,1,...,R\} $ then $m= m+0 \in \{0,1,...,R\}+A $, and if $ m>R $ then $ m-R>0 $ and so $ m-R+k=a \in A$, with $ k \in \{0,1, \ldots,R\} $ i.e. $ 	m=(R-k) +a $ with $ R-k \in \{0,1,...,R\} $. Thus $  \Z_+ = \{0,1,2,..,R\} + A $.
\end{proof}

\begin{corollary} \label{both=}	
	For semicascades both the definitions of almost periodic points coincide.
\end{corollary}
\begin{proof}
	We note that the two conditions in Lemma \ref{equal} are exactly the two kinds of almost	periodicity applied to semicascades.
\end{proof}

We note that the inverse limit construction converts a semicascade to a cascade. This gives some kind of equivalence of minimality in semicascades and cascades.

\begin{theorem}
	
Let $ f:X \to X $ be continuous and surjective. Define $ X^*=: \{(x_0,x_1,\ldots): f(x_i)=x_{i-1}, i \in \N\} $
and let $ f^*:X^* \to X^* $ such that $ f^*(x_0,x_1,\ldots)=(f(x_0),f(x_1),\ldots)=
(f(x_0),x_0,x_1,\ldots) $. Thus $ f^* $ is a homeomorphism giving a cascade $(X^*,f^*)$.

Then $(X, f) $ is minimal if and only if $ (X^*,f^*) $ is minimal.

\end{theorem}

\begin{proof}
If $\pi_0: X^* \to X$ is the projection on the first factor, then $\pi_0(x_0, x_1, \ldots) = x_0$ is a continuous surjection and $\pi_0(f^*(x_0, x_1, \ldots)) = f(x_0) = f(\pi_0(x_0,x_1, \ldots))$ and so $\pi_0$ is a homomorphism. Hence $(X,f)$ is minimal whenever $(X^*,f^*)$ is.

Conversely, let $(X,f)$ be minimal.  We show that $(X^*,f^*)$ is minimal. Let $ x^*=(x_0, x_1,x_2,...) \in X^* $. We show its orbit is dense in $ X^* $. Let
$ {x^*}'=(x_0',x_1',...) \in X^* $. Now let $ m>0  $ and let the sequence $ \{n_i\} $ be such that $ f^{n_i}(x_k) \to x_k' $ for $k = 1, \ldots,m$. Then $ f^{(n_i+j)}(x_k) \to f^j(x_k') $ for $k=  1, \ldots,m) $. That is $ f^n_i(x_l)
\to x_l' $ for $ l=0,1,...m $. Since $ X^* $ has the product topology this proves that $ {x^*}' $ is in the orbit closure of $ x^* $.
	
\end{proof}

\bigskip

\section{Equicontinuity and Distality of Semiflows}

In this section, let $(X,d)$ be a compact metric space.

Recall that the semiflow $(X,T)$ is called surjective if each $t \in T$ is surjective. The TDS $(X,T)$ is said to be \emph{equicontinuous at a point $y\in X$} if for every $\epsilon> 0$ there exists a neighborhood $U$ of $y$ such that $d(tx,ty)<\epsilon$ for every $x\in U$ and every $t\in T$,  and $(X,T)$ is \emph{equicontinuous} if it is equicontinuous for every point of $X$. 

If  $\inf \limits_{t\in T} \ d(tx,ty) = 0$ then the pair $(x,y)$ is called \emph{proximal}. We denote the collection of proximal pairs by $P(X)$, and the TDS $(X,T)$ is distal if and only if $P(X) = \Delta$ - the diagonal in $X \times X$. Thus the TDS $(X,T)$ is  \emph{distal } if for every $(x,y)\in X\times X$ with $x \neq y$ satisfies  $\inf \limits_{t\in T} \ d(tx,ty) > 0$.  Note that if $ (X,T) $ is distal then each $ t \in T $ is injective.

\bigskip

We recall from Corollary \ref{admit} that every semiflow admits semiflow almost periodic points. Almost periodicity has deep connections with distality and proximality. We look into some such connections for semiflows.

\begin{theorem}
	Let $(X,T)$ be a semiflow. If $(x,y) \in X \times X \setminus \Delta$ is semiflow almost periodic for $(X \times X, T)$, then the pair $(x,y)$ is distal. 
\end{theorem}
\begin{proof}
	Let $(x,y) \in X \times X \setminus \Delta$ be semiflow almost periodic and $(x,y) \in P(X)$. Then there is a $z \in X$ and a net $\{t_\alpha\}$ such that $t_\alpha(x,y) \to (z,z)$. Since $\overline{T(x,y)}$ is minimal, $\overline{T(x,y)} = \overline{T(z,z)} \subset \Delta$, giving $x = y$ - a contradiction.
	
	So $(x,y)$ is distal.
\end{proof}

\begin{definition}
	Let $(X,T)$ be a semiflow. Then $A \subset X$ is called \emph{semiflow almost periodic set} if $z \in X^{|A|}$ with $range \ z = A$ where $|A|$ denotes the cardinality of $A$, then $z$ is a semiflow almost periodic point of the semiflow $(X^{|A|}, T)$.
\end{definition}

Note that Zorn's lemma guarantees the existence of  \emph{maximal semiflow almost periodic sets} for any semiflow $(X,T)$.

\begin{theorem}
	Let $(X,T)$ be a semiflow and let $x \in X$. Then there exists a semiflow almost periodic point $x^* \in X$ such that $(x,x^*) \in P(X)$.
\end{theorem}

\begin{proof}
	If $x$ itself is semiflow almost periodic then $x^* = x$. 
	
	We assume that $x$ is not semiflow almost periodic. Let $A \subset X$ be a maximal semiflow almost periodic set. Let $z \in X^{|A|}$ with $range \ z = A$.
	
	Then $(x,z) \in X \times X^{|A|}$. Let $M \subset \overline{T(x,z)}$ be minimal in $X \times X^{|A|}$, and $(x',z') \in M$. Since $z$ is semiflow almost periodic $z' \in \overline{Tz}$ which is minimal, so there exists a net $\{t_\alpha\}$ such that $t_\alpha(z') \to z$. By passing to a subnet if needed, let $t_\alpha(x') \to x^*$ shich is a semiflow almost periodic in $X$. Then $t_\alpha(x',z') \to (x^*,z)$ in $M$. Then $A \cup \{x^*\} = range (x^*,z)$ is a semiflow almost periodic set. By maximality of $A$, it follows that $x^* \in A$.
	
	Now since $(x^*,z) \in \overline{T(x,z)}$, there is a net $\{s_\beta\}$ such that $s_\beta(x,z) \to (x^*,z)$. Now $s_\beta(x) \to x^*$ and $s_\beta(z) \to z$. But $x^* = \pi_{x^*}(z) \in A$ where $\pi_{x^*} $ is the projection on the $x^*$ coordinate in $z$, and hence $s_\beta(x^*) = s_\beta(\pi_{x^*}(z)) = \pi_{x^*}(z) = x^*$.
	
	Thus, $s_\beta(x,x^*) = (x^*,x^*) \in \Delta$ and so $(x,x^*) \in P(X)$.
	
\end{proof}

This also proves a generalization of a classical theorem of Ellis:

\begin{theorem} \cite{tdes} \label{ellis}
	Let $(X,T)$ be a distal semiflow. Then  $ \overline{Tx} $ is minimal for every $x \in X$.
	
	In particular, distal semiflows are  union of minimal subsemiflows.
\end{theorem}

This leads to an interesting consequence.

\bigskip

\begin{proposition} \cite{tdes}
	Suppose (X,T) a distal semiflow. Then  (X,T) is surjective.
\end{proposition}

We give a more concise proof here.

\begin{proof}
	Since (X,T) is distal, for each $t \in T$ the semicascade $(X,t)$ is also distal. Thus every $(X,t)$ is a disjoint union of minimal subsetsand so every $x \in X$ is recurrent. Hence by Lemma \ref{rec}, $t$ is surjective on $X$. So for each $t\in T$, $t$ is surjective on $X$. Hence $(X,T)$ is a surjective semiflow.
\end{proof}

Thus if $ (X,T) $ is distal, every $ t \in T $ is injective and hence invertible and so $T$ becomes a semigroup of homeomorphisms. Thus for $ [T] - $  the group generated by $ T $, the semiflow $(X,T)$ generates the flow $(X, [T])$.

\begin{lemma} \label{minequal}
	Let $(X,T)$ be  distal semiflow and $ [T] $ be the group generated by the semigroup $ T $ of homeomorphisms acting on $X$.
	Then $\overline{[T]x} = \overline{Tx} $, i.e., $\overline{[T]x}$ is also minimal for every $x \in X$.
\end{lemma}

\begin{proof}
	Let $ \overline{Tx} = M$ be minimal in $ X$. Let $y = t^{-1}x \in X $ for some $t \in T$. Then $\overline{Ty} = N$ is also minimal in $X$ and $x = ty \in N$. It follows that the minimal sets $M = N$ which implies that $t^{-1}M \subseteq M$ for every $t \in T$. Thus for each $x \in X $, we have $\overline{Tx} \subseteq \overline{[T]x} \subseteq \overline{Tx}$ is minimal.
\end{proof}

\bigskip

This leads to  a very captivating observation, for which we again give a more concise proof.

\begin{theorem} \cite{tdes} \label{distal}
	For a semiflow $(X,T)$ the following conditions are equivalent:
	\begin{enumerate}
		\item $(X,T)$ is distal.
		\item $(X,[T])$ is distal.
	\end{enumerate}
	
\end{theorem}

\begin{proof}

Since  $ (X ,T) $  is distal, each $t \in T$ is bijective and the product semiflow $(X \times X, T)$ is also distal. By Lemma \ref{minequal} for $(x,y) \in X \times X$ we have $\overline{T(x,y)} \subseteq \overline{[T](x,y)} \subseteq \overline{T(x,y)}$ is minimal. Thus every $(x,y) \in X \times X$ is almost periodic for the flow $(X \times X,[T])$  which immediately implies that $ (X,[T]) $ is distal.

The converse is clear.

\end{proof}

\begin{remark}
	Suppose $ T $ is a group and $ (X,T) $ is distal. Let $ S $ be a subsemigroup of $ T $.
	Then $  (X,S) $ is distal.
	
	We note that the converse holds here too. All distal semigroup
	actions arise this way. 
\end{remark}

\begin{remark}  We note that for a family $\{(X_\lambda,T)\}_{\lambda \in \Lambda}$ of semiflows, the product semiflow $(\prod \limits_{\lambda \in \Lambda} \ X_\lambda, T)$ is distal if and only if every $(X_\lambda,T)$, $\lambda \in \Lambda$, is distal.

Theorem \ref{ellis} further specifies that in such a case the product semiflow $(\prod \limits_{\lambda \in \Lambda} \ X_\lambda, T)$ is pointwise semiflow almost periodic. \end{remark}

It is known that  equicontinuous flows are distal and the converse does not hold.
We   have the following connecting distal and equicontinuous semiflows:

\begin{theorem}  \label{dai}
	Let $(X,T)$ be equicontinuous and surjective semiflow. Then $(X,T)$ is distal.
\end{theorem}

\begin{proof} We note that the equicontinuous and surjective  $T$ is in fact a semigroup or monoid of homeomorphisms on $X$.
	
	We first see that for equicontinuous and surjective $(X,T)$ we have that every $x \in X$ is recurrent for the semicascade $(X,t)$, for all $t \in T$.
	
	Let $x \in X$, $t \in T$ and recursively define $x_i \in X$ such that $ t(x_1) = x $, and $t(x_n) = x_{n-1}$ for all $n > 1$. Let $\delta > 0$ correspond to a given $\epsilon > 0$ in the definition of equicontinuity of $(X,t)$. Pick $m,k >0$ so that $d(x_m, x_{m+k}) < \delta$. Then by equicontinuity, $d(t^{m+k}(x_m), t^{m+k}(x_{m+k})) = d(t^k(x), x) < \epsilon$, proving our assertion.
	
	\bigskip
	
	Now for the surjective, equicontinuous $(X,T)$   let $(x,y) \in P(X)$. Choose $ \epsilon > 0$ such that $d(x,y)/2 > \epsilon$ and let $\delta > 0$ correspond to this $ \epsilon $ in the definition of equicontinuity of $(X,T)$.  Let $t \in T$ be such that $d(t(x), t(y)) < \delta$. Then  $d(t^k(x), t^k(y)) < \epsilon$ for all $k \in \N$. 
	
	The product semicascade $(X \times X, t \times t)$ is also equicontinuous and surjective  and so $(x,y)$ is a recurrent point for this semicascade. 
	
	But  $d(t^k(x), t^k(y)) < \epsilon$ for all $k \in \N$ contradicting the recurrence of $(x,y) $ in $(X \times X, t \times t)$. Hence $x = y$, proving the distality of $(X,T)$.	
\end{proof}

We thank Xiongping Dai for the arguments in the above proof.

\begin{theorem}
	Let $ (X,T) $ be an equicontinuous, surjective semiflow. Then $ (X,[T]) $ is equicontinuous.
	
\end{theorem}

\begin{proof}
	
	The semiflow  $ (X,T) $ is equicontinuous and surjective. Hence, by Theorem \ref{dai} $(X,T)$ is distal, and so is $(X, [T])$ by Theorem \ref{distal}. 
	
	Recall that a flow is equicontinuous if and only if the regionally proximal relation is trivial. So let $ (x,y) \in RP_{[T]} $. That is there are nets $ \{x_n\} $ and $ \{y_n\} $ in $ X $ and $ \{s_n\} $ in $ [T] $ with $ x_n \to x $, $ y_n \to y $ and $ z \in X $ with
	$ s_n(x_n,y_n) \to (z,z) $ for some $  z \in X $. By Lemma \ref{minequal}, we have $ (z,z) \in \overline{[T](x_n,y_n)} = \overline{T(x_n,y_n)} $ so there are $ t_n \in T $ with $ t_n(x_n,y_n) \to (z,z) $.
	Now it follows from the assumed equicontinuity of the $ T $ action that $ t_n(x,x_n) \to (z,z) $ and $ t_n(y,y_n) \to (z,z) $, so we have $ t_n(x,y) \to (z,z) $. That is  $ (x, y) \in P_T(X) $. But as we have noted  $ (X,T) $ is distal 	and so $ x=y $.

\end{proof}

\vspace{12pt}
\bibliography{xbib}

\begin{thebibliography}{99}
	
	
	\bibitem{go} W. H. Gottschalk, Almost Periodic Points with Respect to Transformation Semi-Groups,  Annals of Mathematics ,  47(1946),  762-766.
	
	\bibitem{gh} W. H. Gottschalk and G. A. Hedlund, Topological dynamics, {\it Amer.	Math. Soc. Colloquium Publications}  (1955).
	
	
	\bibitem {tdes}
	 Anima Nagar and Manpreet Singh, Topological Dynamics of  Enveloping Semigroup, {\it  SpringerBriefs in Mathematics}, Springer  (2023).
	
		
	
	
\end{thebibliography}

\end{document}